\documentclass[11pt,a4paper,reqno]{article}
\usepackage{amsmath}
\usepackage{amsthm}
\usepackage{enumerate}
\usepackage{amssymb}

\usepackage{verbatim}
\usepackage{url}
\usepackage[latin1]{inputenc}

\usepackage{caption}
\usepackage{graphicx,color}
\def\Z{{\mathbb Z}}

\def\Bc{{\mathcal B}}

\def\Tc{{\mathcal T}}

\def\be{\begin{equation}}
\def\ee{\end{equation}}
\def\bea{\begin{equation*}}
\def\eea{\end{equation*}}
\def\bal{\begin{aligned}}
\def\eal{\end{aligned}}

\def\eps{\varepsilon}

\def\Pr{{\mathbb P}}
\DeclareMathOperator{\E}{{\mathbb E}}

\newtheorem*{thm}{Theorem}
\newtheorem*{cor}{Corollary}

\theoremstyle{remark}
\newtheorem*{remark}{Remark}

\theoremstyle{definition}

\begin{document}

\title{A temporal perspective on the rate of convergence in first-passage percolation under a moment condition}
\date{\today}
\author{Daniel Ahlberg\thanks{IMPA, Estrada Dona Castorina 110, 22460-320 Rio de Janeiro, Brazil.} \thanks{Department of Mathematics, Uppsala University, SE-75106 Uppsala, Sweden.}}
\maketitle

\begin{abstract}
We study the rate of convergence in the Shape Theorem of first-passage percolation, obtaining the precise asymptotic rate of decay for the probability of linear order deviations under a moment condition. Our results are stated for a given time and complements recent work by the same author, in which the rate of convergence was studied from the standard spatial perspective.
\end{abstract}

\section{Introduction}

Consider first-passage percolation on the $\Z^d$ nearest-neighbour graph for $d\ge2$. Large deviations were first studied in the context of first-passage percolation in the 1980s, in a pioneering work of Grimmett and Kesten~\cite{grikes84}. In this work, together with the subsequent work of Kesten~\cite{kesten86}, the authors investigate the rate of convergence of travel times to the so-called time constant, and provide necessary and sufficient conditions for exponential decay for the probability of linear order deviations. To obtain the exponential decay one requires finite moment of exponential order of the passage times.

It was only recently, in~\cite{ahlberg-8}, that large deviations in the regime of polynomial decay of the probability tails were studied. This regime is highly interesting since it is in this regime that strong laws such as the Subadditive Ergodic Theorem, due to Kingman~\cite{kingman68}, and the Shape Theorem, whose precise conditions were found by Cox and Durrett~\cite{coxdur81}, cease to hold. In this paper we build on this recent work and offer a temporal perspective on the rate of convergence in the Shape Theorem. We complement the results of~\cite{ahlberg-8} with estimates on the probability of linear order deviations at a given time. The results are sharp in the regime of polynomial decay.

We will assume throughout that the edges of the $\Z^d$ lattice are assigned independent non-negative random weights from a distribution satisfying $F(0)<p_c(d)$. Here $p_c(d)$ denotes the critical probability for bond percolation on $\Z^d$. The weights induce a random metric structure to $\Z^d$, and we denote the distance metric by $T(x,y)$ for $x,y\in\Z^d$. Let
$$
\mu(z):=\lim_{n\to\infty}\frac{T(0,nz)}{n},
$$
which is known to exist in probability without further assumption on the edge weights. The limit $\mu(\cdot)$ is usually known as the \emph{time constant}. Richardson~\cite{richardson73} was the first to realize that the above convergence holds in all directions simultaneously, in that $\Bc_t:=\{z\in\Z^d:T(0,z)\le t\}$ obeys a law of large numbers.

Our assumption on the edge weights are known to imply that $\mu(z)\neq0$ for all $z\neq0$ (see~\cite[Theorem~6.1]{kesten86}). Consequently, the deterministic set $\Bc^\mu_t:=\{z\in\Z^d:\mu(z)\le t\}$ is finite for all $t\ge0$. Cox and Durrett~\cite{coxdur81} provided the precise condition for Richardson's result to hold. Let $Y$ denote the minimum of the $2d$ weights assigned to the edges adjacent to the origin. The result of Cox and Durrett says that if $\E[Y^d]<\infty$, then, almost surely, for every $\eps\in(0,1)$ 
\be\label{eq:shapeeq}
\Bc^\mu_{(1-\eps)t}\subset\Bc_t\subset\Bc^\mu_{(1+\eps)t}
\ee
for large enough $t$. The given moment condition is also necessary for this conclusion.

In this paper we aim to investigate in further detail the probability that~\eqref{eq:shapeeq} fails, and that $\Bc_t$ deviates significantly from its asymptotic rate of growth. Let
$$
\Tc_\eps:=\{t\ge0:\text{either inclusion in~\eqref{eq:shapeeq} fails}\}.
$$
We will study the behaviour of $\Pr(t\in\Tc_\eps)$ for large $t$ and fixed $\eps\in(0,1)$. Cox and Durrett's result implies that $\Tc_\eps$ is almost surely bounded if and only if $\E[Y^d]<\infty$. In fact, $\Pr(t\in\Tc_\eps)$ may even be bounded away from zero unless $\E[Y^d]<\infty$, as we remark upon at the very end. Instead we assume that $\E[Y^d]<\infty$ holds, and provide the following results. The first result shows that in the regime of polynomial decay of the tails, $\Pr(t\in\Tc_\eps)$ is governed by the tails of $Y$.

\begin{thm}\label{thm}
Assume that $d\ge2$, $F(0)<p_c(d)$ and $\E[Y^d]<\infty$. Then, for every $\eps>0$ and $q\ge1$, there is $M=M(\eps,d,q)$ such that
$$
M^{-1}\,t^d\,\Pr(Y>t)\,\le\,\Pr(t\in\Tc_\eps)\,\le\, M\,t^d\,\Pr(Y>t/M)+M\,t^{-q}\quad\text{for }t\ge1.
$$
\end{thm}

Using Markov's inequality we find that for each $\alpha>0$, the condition $\E[Y^{d+\alpha}]<\infty$ implies that $\Pr(t\in\Tc_\eps)=O(t^{-\alpha})$.
Another consequence of the theorem is the following characterization.

\begin{cor}
Assume that $F(0)<p_c(d)$. For every $\alpha\ge0$, $\eps>0$ and $d\ge2$,
$$
\E[Y^{d+\alpha}]<\infty\quad\Leftrightarrow\quad\int_0^\infty t^{\alpha-1}\,\Pr(t\in\Tc_\eps)\,dt<\infty.
$$
\end{cor}

Using Fubini's theorem we see that
$$
\E|\Tc_\eps|\,=\,\E\int_0^\infty1_{\{t\in\Tc_\eps\}}\,dt\,=\,\int_0^\infty\Pr(t\in\Tc_\eps)\,dt,
$$
where $|\cdot|$ denotes Lebesgue measure. Hence the conclusion of the corollary was for $\alpha=1$ known already in~\cite[Theorem~2]{ahlberg-8}, but the statement for general $\alpha>0$ was previously unknown. The proofs we present below will be based on results and ideas from~\cite{ahlberg-8}.

\section{Proof}

We first recall some results that will be required for the analysis: For every $\eps>0$ and $d\ge2$ there exists constants $M=M(\eps,d)$ and $\gamma=\gamma(\eps,d)$ such that
\be\label{eq:lowerLDE}
\Pr\big(T(0,z)-\mu(z)<-\eps x\big)\,\le\, M\,e^{-\gamma x}\quad\text{for all }z\in\Z^d\text{ and }x\ge|z|.
\ee
If, in addition, $\E[Y^\alpha]<\infty$ for some $\alpha>0$ and $q\ge1$, then we may choose $M=M(\alpha,\eps,d,q)$ such that
\be\label{eq:upperLDE}
\Pr\big(T(0,z)-\mu(z)>\eps x\big)\,\le\, M\,\Pr(Y>x/M)+M\,t^{-q}\quad\text{for all }z\in\Z^d\text{ and }x\ge|z|.
\ee
The former statement was first proved by Grimmett and Kesten~\cite{grikes84,kesten86} for coordinate directions, and later extended in~\cite{ahlberg-8}. The latter statement is original in~\cite{ahlberg-8}.

Fix $\eps>0$ and $q\ge1$. We wish to estimate the decay of $\Pr(t\in\Tc_\eps)$. Define thus
\bea
\begin{aligned}
A(z)\,&:=\,\{t\ge0: T(0,z)>t\text{ and }\mu(z)\le(1-\eps)t\},\\
B(z)\,&:=\,\{t\ge0: T(0,z)\le t\text{ and }\mu(z)>(1+\eps)t\}.
\end{aligned}
\eea
Note that $t\in\Tc_\eps$ if and only if $t\in A(z)\cup B(z)$ for some $z\in\Z^d$.

\subsection{The lower bound}

We begin with the lower bound. Let $D$ denote the set of all $z\in\Z^d$ such that $\mu(z)\le(1-\eps)t$ and whose $\ell^1$-distance from the origin is even. Let $Y(z)$ denote the minimum weight among the $2d$ edges adjacent to $z$, and note that the $Y(z)$'s are independent for $z\in D$, as points in $D$ are at $\ell^1$-distance at least 2. Clearly $A(z)\subset\Tc_\eps$ for every $z\in\Z^d$. Since also $T(0,z)\ge Y(z)$ we obtain
\be\label{eq:lower}
\begin{aligned}
\Pr(t\in\Tc_\eps)\;&\ge\;\Pr\big(T(0,z)>t\text{ for some }z\in D\big)\\
&\ge\;\Pr\big(Y(z)>t\text{ for some }z\in D\big).
\end{aligned}
\ee

An application of Cachy-Schwarz's inequality shows that any non-negative random variable $X$ satisfies $\Pr(X>0)\ge\E[X]^2/\E[X^2]$, and if $X$ is binomially distributed with parameters $n$ and $p$, then a further lower bound is given by $np/(1+np)$. Applying this to~\eqref{eq:lower} leaves us with
\be\label{eq:lower2}
\Pr(t\in\Tc_\eps)\,\ge\,\frac{|D|\,\Pr(Y>t)}{1+|D|\,\Pr(Y>t)}.
\ee
By assumption $\E[Y^d]<\infty$, which implies that $t^d\,\Pr(Y>t)\le \E[Y^d]$ via Markov's inequality. As the set $D$ grows as $t^d$, we obtain the required lower bound on $\Pr(t\in\Tc_\eps)$.

\subsection{The upper bound}

We now continue with the upper bound. The union bound leaves us with the upper bound
$$
\Pr(t\in\Tc_\eps)\,\le\sum_{\mu(z)\le (1-\eps)t}\Pr\big(T(0,z)>t\big)\,+\!\sum_{\mu(z)>(1+\eps)t}\Pr\big(T(0,z)\le t\big).
$$
We treat the two sums separately. Using~\eqref{eq:lowerLDE} we find $M_1=M_1(\eps,d)$ and $\gamma_1=\gamma_1(\eps,d)$ such that
\bea
\bal
\sum_{\mu(z)>(1+\eps)t}\Pr(T(0,z)\le t)\;&\le\,\sum_{\mu(z)\ge(1+\eps)t}\Pr\left(T(0,z)-\mu(z)\le-\frac{\eps}{1+\eps}\mu(z)\right)\\
&\le\;\sum_{\mu(z)\ge t}M_1\,e^{-\gamma_1 \mu(z)},
\eal
\eea
which is at most $M_2\,e^{-\gamma_2t}$ for $t\ge1$, for some constants $M_2$ and $\gamma_2$. Using~\eqref{eq:upperLDE} we may find $M_3=M_3(\eps,d,q)$ such that
\bea
\bal
\sum_{\mu(z)\le (1-\eps)t}\Pr\big(T(0,z)>t\big)\;&\le\,\sum_{\mu(z)\le(1-\eps)t}\Pr\big(T(0,z)-\mu(z)>\eps t\big)\\
&\le\;\sum_{\mu(z)\le t}M_3\,\Pr(Y\ge t/M_3)+\frac{M_3}{t^{d+q}}.
\eal
\eea
Since the cardinality of the set $\Bc^\mu_t$ grows at the order of $t^d$ we obtain, for some $M_4$,
$$
\Pr(t\in\Tc_\eps)\,\le\, M_4\,t^d\,\Pr(Y\ge t/M_3)+M_4\,t^{-q}+M_2\,e^{\gamma_2t},
$$
for all $t\ge1$, which gives the required upper bound.

\begin{remark}
Finally, we mention that the condition $\E[Y^d]<\infty$ cannot be relaxed in general. Consider for instance the case when $\Pr(Y>t)=t^{-d}$ for $t\ge1$. Then $\E[Y^\alpha]$ is finite for all $\alpha<d$ and infinite for $\alpha=d$, but the bound in~\eqref{eq:lower2} shows that $\Pr(t\in\Tc_\eps)$ is bounded away from zero.
\end{remark}

\newcommand{\noopsort}[1]{}\def\cprime{$'$}

\end{document}